\documentclass[12pt,a4paper]{article}
\usepackage{ifthen} 
\usepackage{calc} 

\usepackage[german,english]{babel} 

\usepackage[all]{xy}
\usepackage{amsmath}
\usepackage{amssymb}
\usepackage{latexsym}
\usepackage{enumerate}
\usepackage{theorem}
\usepackage{mathrsfs}
\usepackage{helvet}

\usepackage{lastpage} 
\usepackage{prettyref} 
\usepackage[pagebackref]{hyperref} 

\input xy

\newboolean{finalversion} 
\setboolean{finalversion}{true}

\theoremheaderfont{\scshape}

\newrefformat{chap}{\hyperref[{#1}]{Chapter~\ref*{#1}}}
\newrefformat{sec}{\hyperref[{#1}]{Section~\ref*{#1}}}

\newtheorem{theorem}{Theorem}[section]
\newrefformat{thm}{\hyperref[{#1}]{Theorem~\ref*{#1}}}

\newrefformat{def}{\hyperref[{#1}]{Definition~\ref*{#1}}}

\newrefformat{lem}{\hyperref[{#1}]{Lemma~\ref*{#1}}}

\newrefformat{prop}{\hyperref[{#1}]{Proposition~\ref*{#1}}}

\newrefformat{cor}{\hyperref[{#1}]{Corollary~\ref*{#1}}}
{\theorembodyfont{\upshape}
\newtheorem{examplecore}[theorem]{Example}}
\newrefformat{ex}{\hyperref[{#1}]{Example~\ref*{#1}}}

\newrefformat{rem}{\hyperref[{#1}]{Remark~\ref*{#1}}}

\newcommand{\mylabel}[1]{\label{#1}\ifthenelse{\boolean{finalversion}}{
  }{\marginpar{\tiny #1}}}  

\title{A Short Note on Polynomial Automorphisms}

\author{Stefan G\"unther}

\begin{document}

\maketitle
\today
\begin{abstract} In this paper, we construct explicitely polynomial automorphisms of affine $n$-space for certain $n\in \mathbb N$\,. More precisely, we construct algebraic subgroups of the general polynomial group $\text{GA}_n(k)$\,, where $k$ is an arbitrary ring of characteristic zero. The formulas are universal and work over any ring of $\mbox{char(0)}$\,.
\end{abstract}
\tableofcontents

\section{Notations and Conventions}
If $k\longrightarrow A$\, is an extension of commutative rings, then $\Omega^{(1)}(A/k)$\, denotes the classical module of Kaehler differentials.\\
 If $N\in \mathbb N$\,, 
$$\mathcal J^N(A/k):=A\otimes_kA/\mathcal I_{\Delta}^{N+1}$$ denotes the $N^{th}$ jet-algebra.\\
 The graded algebra of $N^{th}$-order Kaehler differentials we denote by $\Omega^N(A/k)$\,. If $N=1$\, then 
$$\Omega^1(A/k)=\sum_{n\in \mathbb N_0}\Omega^{(1)}(A/k)^{\otimes sym,n}.$$
The ordinary differential of a function $f\in A$\, we denote by $d^1f$\,.\\
The $n^{th}$ order differential of a function $f\in A$\, we denote by $d^nf$.\\ 
If $(l_{ij})$\, is an $m\times N$-matrix, we let 
$$\lVert l_{ij}\rVert :=\sum_{i,j}l_{ij}\,\,\mbox{and}\,\, \mid l_i\mid:=\sum_{j=1}^Nl_{ij}.$$  
For $k$ a commutative ring, we denote by $\mbox{GA}_n(k)$\, the general polynomial automorphism group of the free $k$-algebra $k[x_1,...,x_n]$\,.\\
By $\mbox{EA}_n(k)$\, we denote the subgroup of elementary polynomial automorphisms and by $T_n(k)$\, the subgroup of triangular polynomial automorphisms, see \cite{PolAut}.\\
\section{Introduction} Although it is known, that the general polynomial group $\text{GA}_n(k)$\, is something like an infinite dimensional algebraic group, besides the linear or affine automorphisms, it is very hard to construct polynomial automorphisms explicitely. There are few examples, besides some obvious constructions such as elements in the elementary polynomial group $EA_n(k)$\,,  mostly in small dimensions, such as the example of Nagata (see \cite{PolAut}).\\
In this paper, we use the higher Kaehler differential algebras $\Omega^N(k[[x_1,...,x_m]]/k)$\, over power series rings in order to construct algebraic subgroups of the general polynomial group $\mbox{GA}_n(k)$\, for $n=m\cdot N$\,. The idea is quite simple and works as follows. For each $N\in \mathbb N$\, the higher Kaehler differential algebras are   an endofunctor 
$$ \Omega^N(-/k):(k-\mbox{Alg})\longrightarrow (k-\mbox{Alg}),$$ where for each $k$-algebra $B$, $\Omega^N(B/k)$\, is actually a $B$-algebra. One can show that the $N^{th}$-Kaehler differential algebra is the solution to a certain universal problem (see \cite{Promotion}).\\In loc. cit., the $N^{th}$  Kaehler differential algebras are constructed on the category of complete $k$-algebras in order to get finitely generated objects. Observe that the classical module $\Omega^{(1)}(k[[x_1,...,x_n]]/k)$\, is no longer finitly generated. One basically wants, as in the $\mathcal C^{\infty}$-case, that for each $p(x_1,...,x_m)\in k[[x_1,...,x_m]]$\,,
$$(*)\,\,d^1p(x_1,...,x_m)=\sum_{i=1}^m\partial^1/\partial^1x_i(p(x_1,...,x_m))\cdot d^1x_i,$$ where we use formal differentiation. In loc. cit., it is proven that one can calculate the $N^{th}$ Kaehler-differential algebras in the complete case as 
\begin{gather}\Omega^N(k[[x_1,...,x_m]]/k)=k[[x_1,...,x_m]][d^jx_i\mid 1\leq j\leq N, i=1,...,m].
\end{gather} There are universal $k$-algebra-homomorphisms
\begin{gather}k[[x_1,...,x_m]]\longrightarrow k[[x_1,...,x_m]][d^jx_i\mid 1\leq j\leq N, 1\leq i\leq m];\\
p(x_1,...,x_m)\mapsto p(x_1,...x_m)+\sum_{j=1}^Nd^jp(x_1,...,x_m).
\end{gather}
Each $k$-algebra automorphism $\phi: A\longrightarrow A$\, induces a $k$-algebra automorphism
\begin{gather}\Omega^N(A/k)\stackrel{\Omega^N(\phi/k),\cong}\longrightarrow \Omega^N(A/k).
\end{gather}
If $(A,\mathfrak{m},k)$\, is a local ring, then the $k$-algebra automorphism  $\Omega^N(\phi/k)$\, on  $\Omega^N(A/k)$\,, keeps the $k$-subalgebra $A$ fixed and restricted to it gives back the automorphism $\phi$\,. The ideal generated by $\mathfrak{m}$\, inside $\Omega^N(A/k)$\, is kept fixed by the automorphism $\Omega^N(\phi/k)$\, and we get induced automorphisms
\begin{gather}
\overline{\Omega^N(\phi/k)}: \Omega^N(A/k)/\mathfrak{m}\cdot \Omega^N(A/k)\stackrel{\cong}\longrightarrow \Omega^N(A/k)/\mathfrak{m}\cdot \Omega^N(A/k).
\end{gather}
 In \cite{Bochner}[Chapter I, Groups of transformations of formal power series] it is proven that one can iterate formal power series over a field of charateristic zero as long as they have no constant term. I.e., given 
 $$f, p_1,p_2,...,p_m\in (x_1,...,x_m)\cdot k[[x_1,...,x_m]],$$ we can form the formal power series 
 $$f(p_1(x_1,...,x_m),p_2(x_1,...,x_m),..., p_m(x_1,...,x_m))\in k[[x_1,...,x_m]].$$
 If 
 $$p_i(x_1,...,x_m)=\sum_{j=1}^ma_{ij}x_j+\,\,\text{higher order terms},$$ then the $k$-algebra homomorphism
 \begin{gather} \phi: k[[x_1,...,x_m]]\longrightarrow k[[x_1,...,x_m]] \\
 f(x_1,...,x_m)\mapsto f(p_1(x_1,...,x_m),...,p_m(x_1,...,x_m))
 \end{gather} is a $k$-algebra automorphism iff 
 $$\mbox{det}(a_{ij})\neq 0.$$
Putting now $A=k[[x_1,...,x_m]]$\, and using the Kaehler differentials  for complete $k$-algebras (see \cite{Promotion}[chapter 6.3, pp.86-96]) together with the explicit formulas (see \cite{Promotion}[chapter 6.4,Theorem 6.55, pp. 96-100]) we obtain automorphisms
\begin{gather}\overline{\Omega^N(\phi/k)}:\\
k[d^1x_i, d^2x_i,...,d^Nx_i\mid 1\leq i\leq m]\stackrel{\cong}\longrightarrow k[d^1x_i,d^2x_i,...,d^Nx_i\mid 1\leq i\leq m].
\end{gather}
Here $d^jx_i$\, is sent to the residue class modulo $\mathfrak{m}$\, of $d^j\phi_i(x_1,...,x_m)$\,. In \cite{Promotion}[chapter 6.4, pp.96-101] an explicit formula is given for the transformation behaviour of the higher differentials $d^j\phi(x_1,...,x_m)$\, generalizing the classical transformation rule for the ordinary differentials (*). If $d^j$\, is given homogenous degree $j$, these transformation formulas are homogenous polynomials of degree $j$. Now, we put $y_{ij}=d^jx_i$\, and $\overline{\Omega^N(\phi/k)}$\, is then a polynomial automorphism of $k[y_{ij}\mid 1\leq i\leq m, 1\leq j\leq N]$\,, i.e., of affine $N\cdot m$-space. Observe that if we put $N=1,$\, we get a linear automorphism. That's why, ordinary Kaehler differentials don't give anything new. \\
Because $\Omega^N(-/k)$\, is a functor, we get a homomorphism of groups
$$\alpha:\mbox{Aut}_k(k[[x_1,...,x_m]]/k)\longrightarrow \mbox{GA}_n(k).$$ In the transformation formulas, that we will give explicitely later, in the $j^{th}$ transformation rule, only partial derivatives of the $\phi_i(\underline{x})$\, up to total order $j$ enter. Because, up to a known scalar factor, these partial derivatives equal the coefficients of the power series $\phi_i,$\, we can truncate our power series above total degree $N$ and get a homomorphism of groups,
$$\alpha_{N,m}:\mbox{Aut}_k(k[[x_1,...,x_m]]/(x_1,...,x_m)^{N+1})\longrightarrow \mbox{GA}_{N\cdot m}(k).$$
\section{Higher Kaehler differentials algebras}
 The construction and the properties of the classical Kaehler differential module $\Omega^{(1)}(A/k)$\, for a ring extension $k\longrightarrow A$\, are well known. It comes along with a $k$-linear derivation $d^1: A\longrightarrow \Omega^{(1)}(A/k)$\, possessing a universal property for general  $k$-linear derivations $t: A\longrightarrow M$\, for an $A$-module $M$. Now, there is the notion of a higher $k$-linear derivation 
 $$t=(t_0,t_1,...,t_n): t_0:A\longrightarrow B,$$ where $B$ is an $A$-algebra, with the property 
 $$t_k(ab)=\sum_{i+j=k}t_i(a)\cdot t_j(b), k=0,...,n$$
  which generalizes the usual Leibnitz rule 
 $$t_1(ab)=a\cdot t_1(b)+b\cdot t_1(a)=t_0(a)\cdot t_1(b)+ t_1(a)\cdot t_0(b).$$
  The fact that $A$-algebras $B$ are used instead of $A$-modules $M$\, takes account for the fact that higher order derivation are higher order approximations for the structure of the ring extension $k\longrightarrow A,$\, where one can say that the classical derivations are linear approximations to the structure of the ring extension $k\longrightarrow B.$\, In order to pass from an $A$-module $M$ to an $A$-algebra $B,$\, put $B=A\oplus M$\,, the usual square zero extension. The higher Kaehler differential algebras $\Omega^N(A/k)$\, which are graded $A$-algebras are the solution to the corresponding universal problem of finding the "absolute" $N$-derivation. Observe that these are commutative algebras and not the alternating algebra on $\Omega^{(1)}(A/k)$\,. The higher Kaehler differential algebras $\Omega^N(A/k)$\, come along with a universal $N$-derivation $(d^0,d^1,...,d^N),$\, where $d^i: A\longrightarrow \Omega^{(i)}(A/k)$\, maps $A$ into the $i^{th}$ graded piece of $\Omega^N(A/k)$\,. The first graded piece is the classical module of Kaehler differentials $\Omega^{(1)}(A/k)$\,. The $k$-linear maps $d^i$\, assemble to a $k$-algebra homomorphism $d: A\longrightarrow \Omega^N(A/k)$\, as follows easily from the properties of higher derivations. For a detailed account of this theory see \cite{Berger} or \cite{Promotion}. In loc.cit., it is shown that for the free polynomial $k$-algebra $A=k[x_1,...,x_m]$\,, the Kaehler differential algebra $\Omega^N(k[x_1,...,x_m]/k)$\, is the free $k[x_1,...,x_m]$-algebra on the free generators $d^jx_i, i=1,...,m, j=1,...,N.$\, As each polynomial $p(x_1,...,x_m)$\, is sent via $d$ to $p+d^1p+d^2p+...+d^Np,$\, it is clear that the higher differentials $d^np(x_1,...,x_m)$\, must satisfy a universal transformation rule in analogy to the classical rule for $d^1p$\,, where universal means that the formula is independent of the ground ring $k$. We have
 $$d^1(p)=\sum_{i=1}^m\partial^1/\partial^1x_i(p)\cdot d^1x_i.$$ 
 If each generator $d^jx_i$\, gets degree $j$\,, the transformation formulas  for $d^nf$\, are homogenous in degree $n$ with coefficients the  higher partial derivatives of $f$ up to total order $n$\,.\\
 The Kaehler differential algebras $\Omega^N(A/k)$\, are finitely generated $A$-algebras if $A$ is a finitely generated $k$-algebra. It is possible to write down a presentation of $\Omega^N(A/k)$\, out of a given presentation of the $k$-algebra $A$. 
 For nonfinitely generated $A$, in general we do not get finitely generated $A$-algebras $\Omega^N(A/k)$\,, even if $A$ is noetherian. In the classical example, the power series ring $A=k[[x_1,...,x_m]]$\,, even the usual Kaehler differentials $\Omega^{(1)}(A/k)$\, are not a finitely generated $A$-module, mainly because the Leibnitz rule does not allow for termwise formal differentiation. \\
 In order to construct well behaved higher differential algebras in the complete noetherian case there are several possiblities. The first is persued in \cite{Promotion} and uses generalized preadic topologies. The second, more simple way is to make the transformation rules (***) (see section 4.1, Explicite formulas, Theorem 4.2) true for formal power series per definitionem. Then, the higher Leibnitz rules for the higher differentials follow from the rules of partial differentiation. Of course, one has to show that the definition is independent of a choosen formal coordinate system.\\ We then really get what we want, namely, for each $N\in \mathbb N$\,,
$$\Omega^N(k[[x_1,...,x_m]]/k)=k[[x_1,...,x_m]][d^jx_i\mid 1\leq i\leq m, 1\leq j\leq N],$$ whith universal derivation 
\begin{gather*}d: k[[x_1,...,x_m]]\longrightarrow k[[x_1,...,x_m]][d^jx_i\mid 1\leq i\leq m, 1\leq j\leq N]\\
p(x_1,...,x_m)\mapsto p(x_1,...,x_m)+d^1p(x_1,...,x_m)+...+d^Np(x_1,...,x_m).
\end{gather*}
The higher Kaehler differential algebras enjoy a couple of properties. They are well behaved under localization, taking tensor products, taking factor rings, etc. see \cite{Berger} or \cite{Promotion}[chapter 6.2, pp.61-86]. In particular, for each morphism of schemes $X\longrightarrow Y,$\, there is an object $\Omega^N(X/Y)$\, which is a sheaf of graded $\mathcal O_X$-algebras. \\
An important fact one should know is that for each $k$-scheme $X$, the space of $N$-truncated arcs $A^N(X/k)$\, is equal to the relative $\mbox{Spec}$\,, $\mbox{Spec}_X\Omega^N(X/k)$\,.\\
The other generalization of the classical Kaehler differentials are the so called jet algebras $\mathcal J^N(A/k)$\, which are more known to the public. For $k\longrightarrow A$\, an extension of commutative rings, they are defined as $\mathcal J^N(A/k):=A\otimes_kA/\mathcal I_{\Delta}^{N+1}$\, together with a derivation 
$$d: A\longrightarrow \mathcal J^N(A/k),\,\, a\mapsto \overline{1\otimes a}.$$ The ideal $\mathcal I_{\Delta}$\, is the kernel of the multiplication map $\mu_A: A\otimes_kA\longrightarrow A$\, and corresponds geometrically to the ideal sheaf of the diagonal. If one defines $d^1a:=\overline{1\otimes a-a\otimes 1}$\, and writes 
\begin{gather*}\mathcal J^N(k[x_1,..,x_m]/k)=k[x_1,...,x_m][d^1x_1,...,d^1x_m]/(d^1x_1,...,d^1x_m)^{N+1},\\
d:k[x_1,...,x_m]\longrightarrow k[x_1,...,x_m][d^1x_1,...,d^1x_m]/(d^1x_1,...,d^1x_m)^{N+1},\\
p(x_1,...,x_m)\mapsto p(x_1,...,x_m)+d^1p(x_1,...,x_m),
\end{gather*}
then the tranformation rule for the differential $d^1p(x_1,...,x_m)$\, is precisely the $N$-truncated Taylor series expansion for the polynomial $p$\,. In \cite{Promotion}, it is shown, that the higher jet-algebras are universal representing solutions to certain functors, the so called filtered derivations. It is again possible to define these objects in the complete noetherian case, see \cite{Promotion}[chapter 6.5, pp. 101-119]. We do not want to go into further detailes because these jet algebras are not as well suited for our purpose of constructing automorphisms of polynomial algebras, mainly because they are truncated objects. 
\section{The result}
Each $k$-algebra automorphism $\phi: A\longrightarrow A$\, induces a $k$-algebra automorphism
\begin{gather}\Omega^N(A/k)\stackrel{\Omega^N(\phi/k),\cong}\longrightarrow \Omega^N(A/k).
\end{gather}
As already said in the introduction, if $(A,\mathfrak{m},k)$\, is a local ring,  the ideal generated by $\mathfrak{m}$\, inside $\Omega^N(A/k)$\, is kept fix by the automorphism $\Omega^N(\phi/k)$\, and we get induced automorphisms
\begin{gather}
\overline{\Omega^N(\phi/k)}: \Omega^N(A/k)/\mathfrak{m}\cdot \Omega^N(A/k)\stackrel{\cong}\longrightarrow \Omega^N(A/k)/\mathfrak{m}\cdot \Omega^N(A/k).
\end{gather}
Putting $A=k[[x_1,...,x_m]]$\, and using the Kaehler differentials  for complete $k$-algebras (see \cite{Promotion}[chapter 6.3, pp.86-96]), together with the explicit formulas (see \cite{Promotion}[chapter 6.4, pp. 96-101]) we obtain automorphisms
\begin{gather}\overline{\Omega^N(\phi/k)}:\\
k[\overline{d^1x_i}, \overline{d^2x_i},...,\overline{d^Nx_i}\mid 1\leq i\leq m]\stackrel{\cong}\longrightarrow k[\overline{d^1x_i},\overline{d^2x_i},...,\overline{d^Nx_i}\mid 1\leq i\leq m].
\end{gather}
Here $d^jx_i$\, is sent to the residue class modulo $\mathfrak{m}$\, of $d^j\phi_i(x_1,...,x_m)$\,. In \cite{Promotion}[chapter 6.4,pp.96-101] an explicit formula is given for the transformation behaviour of the higher differentials $d^j\phi(x_1,...,x_m)$\, generalizing the classical transformation rule for the ordinary differentials. As in the case $n=1$, it  is homogenous, if $d^j(p)$\, gets homogenous degree $j$\,.  Each $d^n\phi(x_1,...,x_m)$\, is sent to an $A$-linear combination of all monomials $\prod_{i,j}(d^jx_i)^{l_{ij}}$\, of total degree $n$, i.e., with $\sum_{i,j}j\cdot l_{ij}=n$\,. The coefficients are higher order partial derivatives of $\phi$\, up to total order $n$\,.\\
The automorphism $\overline{\Omega^N(\phi/k)}$\, sends  the free generator $\overline{d^nx_i}$  to a $\mathbb C$-linear combination of all monomials of total degree $n$, i.e., to a homogenous polynomial in $\overline{d^1x_i},...,\overline{d^Nx_m}$\, ot degree $n$\,.\\
 Now, we put $y_{ij}=d^jx_i$\, and $\overline{\Omega^N(\phi/k)}$\, is then a polynomial automorphism of $k[y_{ij}\mid 1\leq i\leq m, 1\leq j\leq N]$\,, i.e., of affine $N\cdot m$-space. \\
Because $\Omega^N(-/k)$\, is a functor, we get a homomorphism of groups
$$\alpha:\mbox{Aut}_k(k[[x_1,...,x_m]]/k)\longrightarrow \mbox{GA}_{N\cdot m}(k).$$ In the transformation formulas,  in the $j^{th}$ transformation rule, only partial derivatives of the $\phi_i(\underline{x})$\, up to total order $j$ enter. Because, up to a known scalar factor, these partial derivatives equal the coefficients of the power series $\phi_i,$\, we can truncate our power series above total degree $N$ and get a homomorphism of groups,
$$\alpha_N:\mbox{Aut}_k(k[[x_1,...,x_m]]/(x_1,...,x_m)^{N+1})\longrightarrow \mbox{GA}_{N\cdot m}(k).$$
In \cite{Promotion}, the linear group $\mbox{Aut}_k(k[[x_1,...,x_m]]/(x_1,...,x_m)^{N+1})$\, is called the $N^{th}$ Jacobian linear group because the higher order coefficients of such an automorphism form something like a  higher Jacobian matrix.
In \cite{Promotion}[chapter 9, Proposition 9.2, Proposition 9.4, pp. 148-160] it is proven, using the same idea, that the map
\begin{gather*}\alpha_{(N),m}:\mbox{Aut}_k(k[[x_1,...,x_m]]/(x_1,...,x_n)^{N+1})\longrightarrow\\
 \mbox{Aut}_k(\Omega^{(N)}(k[[x_1,...,x_m]]/k)/(x_1,...,x_m)\cdot \Omega^{(N)}(k[[x_1,...,x_m]]/k))
 \end{gather*}
is an injective algebraic representation of the linear algebraic group 
$$\mathfrak{Jac}_{N,m}(k):=\mbox{Aut}_k(k[[x_1,...,x_m]]/(x_1,...,x_m)^{N+1}).$$
 The right hand side is just the general linear group of the $k$-vector space 
$$W_{N,m}:=\Omega^{(N)}(k[[x_1,...,x_m]]/k)/(x_1,...,x_m)\cdot \Omega^{(N)}(k[[x_1,...,x_m]]/k).$$
Recall that $\Omega^{(N)}(k[[x_1,...,x_m]]/k)$\, is the $N^{th}$ graded piece of the graded $k[[x_1,...,x_m]]$-algebra $\Omega^N(k[[x_1,...,x_m]]/k)$\, and as such a finitely generated free $k[[x_1,...,x_m]]$-module and thus $W_{N,m}$\, is a finite dimensional $k$-vector space. The representation $\alpha_{N,m}$\, is the direct sum over all representations $\alpha_{(n),m}$\,. Thus, we have proved 
\begin{theorem}\mylabel{thm:T1} (Main result)\\
 For each $m,N\in\mathbb N$\, there is an injective homomorphism of groups
$$\alpha_{N,m}: \mbox{Aut}_k(k[[x_1,...,x_m]]/(x_1,...,x_m)^{N+1})\hookrightarrow \mbox{GA}_{N\cdot m}(k).$$
If 
$$\phi=(\phi_1,...,\phi_m)\in \mbox{Aut}_k(k[[x_1,...,x_m]]/(x_1,..,x_m)^{N+1}),$$
 where $\phi_i=\phi_i(x_1,...,x_m)$\, are $N$-truncated power series and $y_{ij}, 1\leq i\leq m, 1\leq j\leq N$\, are the coordinates of affine $N\cdot m$-space, then $\alpha_{N,m}(\phi)(y_{ij})$\, are polynomials of degree at most $N$. Its coefficients are given by universal explicit easy formulas in the coefficients of the power series $\phi_i$\,.
\end{theorem}
\subsection{Explicit formulas}
 In \cite{Promotion}[chapter 6.4, pp.96-98],basically the following theorem is proved.
 \begin{theorem}\mylabel{thm:2} Let $k$ be any base ring of characteristic zero. Then for each $N\in \mathbb N$\,, the $N$-truncated graded ring of complete absolute Kaehler differentials $\Omega^N(k[[x_1,...,x_m]]/k)$\, is isomorphic to the free $k[[x_1,...,x_m]]$-algebra over the indeterminates $d^jx_i, i=1,...,m, j=1,...,N$\,.\\
 The universal $N^{th}$ order derivation is a $k$-algebra homomorphism
 $$k[[x_1,...,x_m]]\longrightarrow k[[x_1,...,x_m]][d^jx_i\mid 1\leq j\leq N, 1\leq i\leq m],$$ where a power series $f=f(x_1,...,x_m)$\, is sent to $f+d^1f+d^2f+...+d^Nf$\,.
 If $k$ is a ring of characteristic zero, then for each $n=1,...,N$\, and each power series 
 $f=f(x_1,...,x_m),$\, we have
 \begin{gather*}(***)\,\,d^n(f)=\sum_{(l_{ij})\in \mathbb N_0^{n\cdot m}, \sum_{i,j}jl_{ij}=n}\frac{1}{\prod_{i,j}l_{ij}!}\cdot\\
  \partial^{\lVert l_{ij} \rVert}/\partial^{\mid l_1\mid}x_1\partial^{\mid l_2\mid}x_2...\partial^{\mid l_m\mid}x_m(f(x_1,...,x_m))\cdot \prod_{i,j}(d^jx_i)^{l_{ij}}.
 \end{gather*}
 \end{theorem}
 To explain something about the formula, the elements $\prod_{i,j}(d^jx_i)^{l_{ij}}$\, with $\sum_{i,j}j\cdot l_{ij}=n$\,  form a free $k[[x_1,...,x_m]]$-basis of the homogenous degree $n$ part $\Omega^{(n)}(k[[x_1,...,x_m]]/k)$ of $\Omega^N(k[[x_1,...,x_m]]/k)$\, and the formula (***) then writes $d^nf(x_1,...,x_m)$\, as a $k[[x_1,...,x_m]]$-linear combination in theses basis elements. Everything is in analogy to the case $n=1$\,. Here, higher partial derivatives of $f$ enter into the formula, namely of total order less then or equal to $n$. The expression $\lVert l_{ij}\rVert$\, denotes the total sum of all $l_{i,j}$\, which can then take each value from $1$ to $n$ (put all $l_{ij}=0$\, but $l_{1n}=1$\, then  
 $$\sum_{i,j}j\cdot l_{ij}=n\,\, \mbox{and}\,\, \lVert l_{ij}\rVert =1,$$
  or all $l_{ij}=0$\,, but $l_{11}=n$\,,and then 
  $$\sum_{i,j}j\cdot l_{ij}=n\,, \mbox{and} \,\,\lVert l_{ij}\rVert =n.$$
 The partial derivatives are taken in a formal sense.\\
  To get an intuitive feeling, we write down the formula $(***)$\, in the special cases $n=m=2$\, and $n=3, m=1$\,. As the interested reader may verify, for $n=1$\,, the formula reduces to the classical formula (*) for the total differential $d^1(f)$\,. For $f(x_1,x_2)\in k[[x_1,x_2]]$\, we have the generalized differentials 
 $$d^2x_1,d^2x_2, (d^1x_1)^2, (d^1x_2)^2\,,\,\mbox{and}\,\,d^1x_1\cdot d^1x_2.$$
 We get 
 \begin{gather*}d^2f(x_1,x_2)= \partial^1/\partial^1x_1(f(x_1,x_2))\cdot d^2x_1+\partial^1/\partial^1x_2(f(x_1,x_2))\cdot d^2x_2+ \\
 \partial^2/\partial^1x_1\partial^1x_2(f(x_1,x_2))\cdot d^1x_1\cdot d^1x_2+ \frac{1}{2}\partial^2/(\partial^1x_1)^2(f(x_1,x_2))\cdot (d^1x_1)^2+\\
 \frac{1}{2}\partial^2/(\partial^1x_2)^2(f(x_1,x_2))\cdot (d^1x_2)^2.
 \end{gather*}
 For $n=3, m=1$\, and $f=f(x)\in k[[x]]$\, we get
 \begin{gather}d^3f(x) =\\
 d^1/d^1x(f(x))\cdot d^3x+ \frac{1}{6}d^3/(d^1x)^3(f(x))\cdot (d^1x)^3+ d^2/(d^1x)^2(f(x))\cdot d^1x\cdot d^2x.
 \end{gather}
 In order to describe our automorphism $\alpha_{N,m}(\phi)$\,, remember, that we have taken the reduction at the maximal ideal $(x_1,...,x_m)$. In our formula $(***)$\, this corresponds to putting all $x_i$\, to zero and considering in the formal power series 
 $$\frac{1}{\prod_{i,j}l_{ij}!}\partial^{\lVert l_{ij} \rVert}/\partial^{\mid l_1\mid}x_1\partial^{\mid l_2\mid}x_2...\partial^{\mid l_m\mid}x_m(f(x_1,...,x_m))$$ 
 only the constant term. It obviously only depends on the coefficients of the power series $f=f(x_1,...,x_m)$\,. We want to explicitely write this down.\\
 Let 
 $$\phi_{\alpha}(x_1,...,x_m)=\sum_{\underline{n}\in \mathbb N_0^m}a_{\underline{n}}^{\alpha}\cdot \underline{x}^{\underline{n}},\,\, \alpha =1,...,m,$$ where we use multiindex notation 
 $$\underline{n}=(n_1,...,n_m),\,\, \underline{x}^{\underline{n}}=x^{n_1}_1\cdot x_2^{n_2}\cdot ...\cdot x_m^{n_m}.$$ 
 The least coefficient of 
 $$\frac{1}{\prod_{i,j}l_{ij}!}\partial^{\lVert l_{ij} \rVert}/\partial^{\mid l_1\mid}x_1\partial^{\mid l_2\mid}x_2...\partial^{\mid l_m\mid}x_m(\phi_{\alpha}(x_1,...,x_m))$$
  is precisely 
  $$c=\frac{\mid l_1\mid !\cdot \mid l_2\mid !\cdot ...\cdot \mid l_m\mid !}{\prod_{i,j}l_{ij}!}\cdot a^{\alpha}_{\underline{\mid l\mid}}.$$
   We can now  write down the universal polynomials for the automorphism
   $\alpha_{N,m}(\phi)$\,.\\ We have 
   \begin{gather}\alpha_{N,m}(\phi)(y_{rs})= \overline{d^s(\phi_r(x_1,...,x_m))}=\\
   \sum_{(l_{ij})\in \mathbb N_0^{s\cdot m}, \sum_{i,j}j\cdot l_{ij}=s}\frac{\prod_{k=1}^m\mid l_k\mid !}{\prod_{i,j} l_{ij}!}\cdot a^r_{\underline{\mid l\mid}}\cdot \prod_{i,j}y_{ij}^{l_{ij}}
   \end{gather}
   We calculate now the  examples $n=m=2$\, and $n=3, m=1$\,. In this case, for $\phi\in \mbox{Aut}_k(k[[x_1,x_2]]/(x_1,x_2)^3)$\,, $\alpha_{2,2}(\phi)$\, is a polynomial automorphism of affine $4$-space $\mathbb A^4_k$\, over $k$\,. As variables we take $y_{11},y_{12},y_{21},y_{22}$\, where $y_{ij}$\, corresponds to the residue of $d^jx_i$\,. Here, $\alpha_{2,2}(\phi)$\, is given by at most quadratic polynomials. We have 
   \begin{gather}\alpha_{2,2}(\phi)(y_{11})=\overline{d^1\phi_1(x_1,x_2)}=\\
   =a^1_{1,0}\cdot y_{1,1}+a^1_{0,1}\cdot y_{2,1}\\
   \alpha_{2,2}(\phi)(y_{21})=a^2_{1,0}\cdot y_{1,1}+a^2_{0,1}\cdot y_{2,1}\\
   \alpha_{2,2}(\phi)(y_{12})=a^1_{1,0}\cdot y_{1,2}+a^1_{0,1}\cdot y_{2,2}+a^1_{2,0}\cdot y_{1,1}^2+ a^1_{0,2}\cdot y_{1,2}^2+ a^1_{1,1}\cdot y_{1,1}\cdot y_{2,1}\\
   \alpha_{2,2}(\phi)(y_{2,2}) =a^2_{1,0}\cdot y_{1,2}+a^2_{0,1}\cdot y_{2,2}+a^2_{2,0}\cdot y_{1,1}^2+\cdot a^2_{0,2}\cdot y_{1,2}^2+ a^2_{1,1}\cdot y_{1,1}\cdot y_{2,1},
   \end{gather}
   where $\phi: k[[x_1,x_2]\stackrel{\cong}\longrightarrow k[[x_1,x_2]]$\,
   is given by 
   \begin{gather}\phi(x_1)=\phi_1(x_1,x_2)=a^1_{1,0}x_1+a^1_{0,1}x_2+a^1_{2,0}x_1^2+a^1_{0,2}x_2^2+a^1_{1,1}x_1x_2;\\
   \phi(x_2)=\phi_2(x_1,x_2)=a^2_{1,0}x_1+a^2_{0,1}x_2+a^2_{2,0}x_1^2+a^2_{0,2}x_2^2+a^2_{1,1}x_1x_2.
   \end{gather}
   The only restriction on the free variables $a^r_{i,j}\in k$\, is 
   $$\mbox{det}\lvert \begin{array}{cc}a_{1,0}^1 a^1_{0,1}\\
                        a^2_{1,0} a^2_{0,1}\end{array}\rvert\neq 0$$ in order to guarantie that $\phi$\, is an automorphism. Here, in this easy example, the scalar factors 
                        $$\frac{\prod_{k=1}^m\mid l_k\mid!}{\prod_{i,j}l_{ij}!}$$ are all equal to one as one easily veryfies.\\ 
   We finally calculate the case $n=3,m=1,$\, where we put $y_j=\overline{d^jx}$\,. In this case, we get a polynomial automorphism of affine $3$-space $\mathbb A^3_k$\, where $\alpha_{3,1}(\phi)$\, is given by at most cubic polynomials.
   We have 
   \begin{gather}\alpha_{3,1}(\phi)(y_1)=a_1\cdot y_1;\\
   \alpha_{3,1}(\phi)(y_2)= a_1\cdot y_2+a_2\cdot y_1^2;\\
   \alpha_{3,1}(\phi)(y_3)=a_1\cdot y_3+a_3\cdot y_1^3+ 2\cdot a_2\cdot y_1\cdot y_2,
   \end{gather}
   where $\phi(x)=a_1x+a_2x^2+a_3x^3\in \mbox{Aut}_k(k[[x]]/(x)^4)$\, and $a_1\neq 0$\,.\\
    One can read off that in the last case we get an automorphism in the triangular polynomial group $T_3(k)$\,.
    For general $n\in \mathbb N,$\, $m=1$\,, it is difficult to write down the formulas for $\alpha_{n,1}(\phi)$\, by hand. Namely, we have to compute all partitions (with repetitions)
    $$(\underbrace{n_1,...,n_1}_{k_1},\underbrace{n_2,...,n_2}_{k_2},...,\underbrace{n_l,...,n_l}_{k_l}),$$ of $n\in \mathbb N$\,, $n=\sum_{i=1}^lk_i\cdot n_i$\, in order to calculate the coefficient in front of 
    $$(d^{n_1}y)^{k_1}\cdot (d^{n_2}y)^{k_2}\cdot ...\cdot (d^{n_l}y)^{k_2}.$$ But as in the case $n=3$\,, one always gets elements in the triangular polynomial group $T_n(k)$\,.\\
    This is a general phenomenon. In \cite{Promotion}[chapter 8.2, pp.136-148], it is shown, that there is a functorial filtration on the higher Kaehler differentials $\Omega^N(A/k),$\,  where functorial means that it is preserved by any $k$-algebra endomorphism $\phi: A\longrightarrow A$\,. In very elementary terms this means that in the transformation rule (***) for $d^nf(x_1,...,x_m)$\, only differentials $d^jx_i$\, with $j\leq n$\, can occur. In our case, where we consider polynomial automorphisms of $\mathbb A^{Nm}_k$\,, if we introduce $m\times m$-blocks, then we can consider the corresponding block-triangular subgroup $TB^m_{Nm}(k)$\, and all our automorphisms will lie in $TB^m_{Nm}(k)$\,. This precisely means that if we consider the original coordinates 
    $$\overline{d^jx_i},\,\, 1\leq i\leq m;\,\, 1\leq j\leq N,$$
     then the affine subspaces 
    $$k[\overline{d^jx_i}\mid 1\leq i\leq m; 1\leq j\leq n_0],\,\, 1\leq n_0\leq N$$
     are kept fixed by our automorphisms. \\
     In order to get more special automorphisms of $\mathbb A^{Nm}_k$\, we can specialize our automorphisms in $\mbox{Aut}_k(k[[x_1,...,x_m]]/(x_1,...,x_m)^{N+1})$\, where only the case $m>1$\, is interesting.  Namely, we may introduce the formal elementary subgroup $FE_m(k)\subset \mbox{Aut}_k(k[[x_1,...,x_m]])$\, (and its $N$-truncated versions), where the last group could be called the formal general group $FG_n(k)$\,. Here, the  formal elementary group  $FE_m(k)$\, is generated by all $\phi=(\phi_1,\phi_2,...,\phi_m)$\, where 
     $$\phi_j(x_1,...,x_m)=x_j\,\forall j\neq i_0\,\, \text{and}\,\, \phi_{i_0}(x_1,...,x_m)=x_{i_0}+\psi(x_1,...,\hat{x_{i_0}},...x_m),$$
     for some $1\leq i_0\leq m$,\,
      where $\psi\in k[[x_1,...,\hat{x_{i_0}},...,x_m]]$\,.
     In this case, we get a "block-elementary" automorphism $\alpha_{n,m}(\phi)$\, namely, the only nontrivial polynomials  are
     $$\overline{d^jx_{i_0}}\mapsto \psi_j(d^rx_i\mid 1\leq r\leq j, 1\leq i\leq m),\,\,1\leq j\leq n$$
     for some polynomials $\psi_j$\, that derive from the polynomial $\psi$\,.\\
    Let the formal triangular group $FT_m(k)\subset FG_m(k)$\, consist of all 
    $$\phi=(\phi_1,...,\phi_m)\,\, \text{such that}\,\, \phi_i=\phi_i(x_1,...,x_i).$$ Then, its image via the homomorphism $\alpha$\, in $GA_{Nm}(k)$\,  clearly is in $T_{Nm}(k)$\,.\\
    Finally, the general linear subgroup $GL_m(k)$\, may be considered as a subgroup of $FG_m(k)$\, by homogenous linear automorphisms. Looking at the transformation formula $(***)$\,, we see that all coefficients in front of mixed monomials get zero and we get a linear automorphism in $GA_{mn}(k)$\,. We have thus constructed a homomorphism of groups 
    $GL_m(k)\hookrightarrow GA_{mn}(k)$\,.   
\section{Generalizations and comments} From the construction principle, it should be clear that working with classical Kaehler differentials only gives linear automorphisms, i.e., nothing new. Our construction gives homogenous polynomial automorphisms, if the coordinate $\overline{d^jx_i}$\, gets degree $j$.  If one defines on $\mbox{GA}_n(k)$\, the structure of an infinite dimensional algebraic group in a resonable way, then our representation
$$\mbox{Aut}_k(k[[x_1,...,x_m]]/(x_1,...,x_m)^{N+1})\hookrightarrow \mbox{GA}_{N\cdot m}(k)$$ should be an algebraic representation. \\
It directly follows from the functoriality of the higher Kaehler differential algebras and their freeness in the case of the polynomial ring, that we get injective homomorphisms of groups 
$$GA_m(k)\hookrightarrow GA_{m\cdot N+m}(k),$$ where 
$\phi\in \mbox{Aut}_k(k[x_1,...,x_m])$\, is simply sent to 
$$\Omega^N(\phi/k)\in \mbox{Aut}_k(k[x_1,...,x_m, d^jx_i, 1\leq i\leq m, 1\leq j\leq N]).$$
Here, also the case $N=1$\, gives something nontrivial, namely a homomorphism of groups
$$GA_m(k)\hookrightarrow GA_{2m}(k).$$ The $(n+i)^{th}$\, variable $d^1x_i$\, is sent to $d^1\phi_i(x_1,...,x_m)$\, if $\phi=(\phi_1,...,\phi_m)$\, and the right hand side of the classical transformation rule 
$$d^1\phi_i(x_1,...,x_m)=\sum_{j=1}^m\partial^1/\partial x_j(\phi_i(x_1,...,x_m))d^1x_j$$ is now to be read as a polynomial in the variables $x_1,...,x_m,d^1x_1,...,d^1x_m$\,.\\
\\
Next, the higher jet bundles or -algebras $\mathcal J^N(A/k)$\, for a ring extension $k\longrightarrow A$, which are more known to the public are not well suited for the construction of polynomial automorphisms. Namely, they are truncated polynomial algebras in case $A=k[x_1,...,x_m]$\, or $A=k[[x_1,...,x_m]]$\,,
$$\mathcal J^N(k[[x_1,...,x_m]]/k)=k[[x_1,...,x_m]][d^1x_1,...,d^1x_m]/(d^1x_1,...,d^1x_m)^{N+1},$$
and, in the untruncated case,
$$\mathcal J^{\infty}(k[[x_1,...,x_m]]/k)=k[[x_1,...,x_m]][[d^1x_1,...,d^1x_m]].$$
 So, reduction modulo the ideal $(x_1,...,x_m)$\, gives in the $N$-truncated case only an automorphism of the Artin-$k$-algebra $k[d^1x_1,...,d^1x_m]/(d^1x_1,...,d^1x_m)^{N+1}$\, and, in the untruncated case an automorphism of the formal power series ring $k[[d^1x_1,...,d^1x_m]],$\, which is nothing but the original automorphism of $k[[x_1,...,x_m]]$\, now written in the variables $d^1x_i$\,. So at least to myself it is not clear how to get polynomial automorphisms out of the mechanism of higher jet algebras.\\
 There are more general higher Kaehler differential algebras. Namely, as explained in \cite{Promotion}, one can replace the grading monoid $\mathbb N$\, by an arbitrary monoid $M$ and the subset $\{1,...,N\}\subset \mathbb N$\, by certain truncation sets $N\subset M$, namely subsets $N\subset M$\, such that for all $n\in N, r+s=n$\, also $r,s\in N$\,. Then, for a ring extension $k\longrightarrow A,$\, there are the  Kaehler diffential algebras $\Omega^M(A/k)$\, and  $\Omega^N(A/k)$\,.  They also can be defined in the complete case and one should have 
 $$\Omega^N(k[[x_1,...,x_m]]/k)=k[[x_1,...,x_m]][d^nx_i\mid n\in N, i=1,...,m].$$
 As an example, take $M=\mathbb N_0^k$,\, $N=\{\underline{n}\mid \sum_{i=1}^kn_i\leq n\}$\,.\\
 Although not explicitely proved, there should also hold a transformation formula (***), namely for all $n\in N,$\,
 \begin{gather}
 d^np(x_1,...,x_m)=\\
 \sum_{(l_{ij})\in \mathbb N_0^{N\cdot m}, \sum_{i,j}jl_{ij}=n}\frac{1}{\prod_{i,j}l_{ij}!}\partial^{\lVert l_{ij}\rVert }/\partial^{\lvert l_1\rvert}x_1...\partial^{\lvert l_m\rvert}x_m (p(x_1,...,x_m))\cdot \prod_{i,j}(d^jx_i)^{l_{ij}},
 \end{gather}
 where, to explain the formula, the sum $\sum_{i,j}j\cdot l_{ij}$\, is taken in the monoid $M$. Of course, we can multiply monoid elements by natural numbers.\\
 Exploiting this example in the same way, we even get non block-triangular polynomial automorphisms.

\bibliography{ShortNote}
\bibliographystyle{plain}
\noindent
\emph{E-Mail-adress:}verb!stef.guenther2@vodafone.de!
\end{document}